\documentclass[reqno]{amsart}

\usepackage{general}
\usepackage{strings}

\numberwithin{equation}{section} \numberwithin{theorem}{section}

\begin{document}

\title{Canonical self-affine tilings by iterated function systems.}
\author{Erin P. J. Pearse}
\address{\scriptsize Department of Mathematics, 310 Malott Hall, Cornell
University, Ithaca, NY 14850-4201} \email{erin@math.cornell.edu}
\thanks{}

\date{\today}



\begin{abstract}
  An iterated function system $\Phi$ consisting of contractive affine
  mappings has a unique attractor $F \subseteq \mathbb{R}^d$ which
  is invariant under the action of the system, as was shown by
  Hutchinson \cite{Hut}.
  This paper shows how the action of the function system naturally
  produces a tiling $\mathcal{T}$ of the convex hull of the attractor.
  These tiles form a collection of sets whose geometry is typically
  much simpler than that of $F$, yet retains key information about
  both $F$ and $\Phi$. In particular, the tiles encode all the
  scaling data of $\Phi$. We give the construction, along with some
  examples and applications.
  The tiling $\mathcal{T}$ is the foundation for the higher-dimensional
  extension of the theory of \emph{complex dimensions} which was
  developed in \cite{FGCD} for the case $d=1$.
\end{abstract}

\maketitle

\section{Introduction}
\label{sec:Introduction}

This paper presents the construction of a self-affine tiling which
is canonically associated to a given self-affine system \simt (as
defined in Def.~\ref{def:self-affine_system}). The term
``self-affine tiling'' is used here in a sense quite different from
the one often encountered in the literature. In particular, the
region being tiled is the complement of the self-affine set \attr
within its convex hull, rather than all of \bRd. Moreover, the tiles
themselves are neither self-affine nor are they all of the same
size; in fact, the tiles may even be simple polyhedra. However, the
name ``self-affine tiling'' is appropriate \bc we will have a tiling
of the convex hull: the union of the closures of the tiles is the
entire convex hull, and the interiors of the tiles intersect neither
each other, nor the attractor \attr. While the tiles themselves are
not self-affine, the overall structure of the tiling
is.\footnote{Technically, the tiling is \emph{subselfaffine} in that
$\simt(\tiling) \ci \tiling$, rather than $\simt(\tiling) =
\tiling$.}

The construction of the tiling begins with the definition of the
generators, a collection of open sets obtained from the convex hull
of \attr. The rest of the tiles will be seen to be images of these
generators under the action of the original self-affine system.
Thus, the tiling \tiling essentially arises as a spray on the
generators, in the sense of \cite{LaPo} and \cite{FGCD}. The tiles
thus obtained form a collection of sets whose geometry is typically
much simpler than that of \attr, yet retains key information about
both \attr and \simt. In particular, the tiles encode all the
scaling data of \simt.

Section \ref{sec:Applications} describes the intended context of the
present results, and indicates how they may be used to develop a
tube formula for the tiling. Section \ref{sec:The Self-Similar
Tiling} gives the tiling construction. Section \ref{ssec:Examples}
illustrates the method with several familiar examples, including the
Koch snowflake curve, Sierpinski gasket and pentagasket. Section
\ref{sec:Properties of the tiling} describes the basic properties of
the tiling and contains the main results of the paper.

\section{Motivation and primary application}
\label{sec:Applications}

The motivation behind the self-affine tiling is to find a means of
extending the work of \cite{FGCD} to higher dimensions. Preliminary
investigations for this project began with \cite{KTF}, and the
present paper, together with \cite{TFCD}, has made significant
inroads. An outline of these results is given in the survey paper \cite{FractalTubeSurvey}. The self-affine tiling provides a natural way to define the
\emph{geometric zeta function} of a self-affine subset of general
Euclidean space \bRd, and thus obtain the \emph{complex dimensions}
of such a set. These terms require some discussion.

The research monograph \mbox{\cite{FGCD}} is an investigation of the
theory of fractal subsets of \bR. The complement of a fractal within
the interval containing it is called a \emph{fractal string} and may
be represented by a sequence of bounded open intervals $L =
\{L_n\}_{n=1}^\iy$, where the interval $L_n$ has length $\ell_n$.
For technical reason, the fractal string may be considered as
\begin{equation} \label{notation:fractal-string}
  \L := \{\ell_n\}_{n=1}^\iy, \q \text{with } \sum_{n=1}^\iy \ell_n < \iy.
\end{equation}
The authors are able to relate geometric and spectral properties of
such objects through the use of zeta functions which encode this
data. The most important such function is the geometric zeta
function \(\gzL(s) = \sum_{n=1}^\iy \ell_n^s\), and the complex
dimensions are defined to be the poles of this function. Prior to
the construction of the self-affine tiling, it was not known how to
define these objects in the higher-\dimnl case, that is, to fractal
subsets of \bRd (for $d>1$).

The poles of the geometric zeta function contain much geometric
information about the underlying fractal, including the dimension
and measurability of the fractal under consideration. Another main
result of \cite{FGCD} is an explicit tube formula, that is, a
formula for the volume $V_\L(\ge)$ of the \ge-\nbd of a fractal
string \L. This formula is roughly just the sum of the residues of
the geometric zeta function at each of the complex dimensions.

The basis for extending these (and other) results to higher
dimensional self-affine sets is a suitable higher-dimensional
analogue of fractal strings: the self-affine tiling developed in
this paper. In \cite{TFCD}, we show how a self-affine tiling \tiling
allows one to define a geometric zeta function for self-affine
subsets of \bRd. Further, we compute an explicit inner tube formula for
$V_\sT(\ge)$ analogous to \cite[Thm. 8.1]{FGCD}, using tools from
geometric measure theory and convexity theory. That is, for $A \ci
\bRd$, we obtain an explicit expression for
\begin{equation}
  V_{A}(\ge) := \vol[d]\{x \in A \suth d(x,\del A)< \ge\},
\end{equation}
where $\vol[d]$ is $d$-\dimnl Lebesgue measure. Each tile in the
self-affine tiling contributes much data to the final formula,
including curvature information for each (topological) dimension
$0,1,\dots,d$, of the tile. Thus, the present construction is used
in an essential manner in \cite{TFCD} to obtain the following result.

\begin{theorem}
  \label{thm:tiling_tube_formula}
  The $d$-\dimnl volume of the inner tubular \nbd of a tiling \tiling
  is given by the following distributional explicit formula:
  \begin{align}
    \label{eqn:preview of main result}
    V_\tiling(\ge)
    &= \sum_{\gw \in \DT} \negsp[1] \res{\gzT(\ge,s)}
     = \sum_{\gw \in \DT} \negsp[1] c_\gw \ge^{d-\gw}.
  \end{align}
\end{theorem}

In this formula, \gzT is the geometric zeta \fn of the self-affine
tiling whose residues define the \coeffs $c_\gw$. The \mero
\dist-valued \fn \gzT is defined by the combinatorics of the scaling
ratios of \simt, and various geometric properties of the generators.
In particular, \gzT incorporates the $0$-\dimnl through $d$-\dimnl
curvatures of the generators of the self-affine tiling. The
sum in \eqref{eqn:preview of main result} is taken over the set of
complex dimensions $\DT$, that is, the set of poles of \gzT. Further
discussion of these topics, however, is beyond the scope of the
current paper; please see \cite{TFCD}.

Additionally, under certain \conds, Theorem~\ref{thm:tiling_tube_formula} can also be used to obtain the volume of the exterior tubular \nbd of a fractal set, that is, an explicit expression for
\begin{equation}
  V_{\attr}(\ge) := \vol[d]\{x \notin \attr \suth d(x,\del \attr)< \ge\}.
\end{equation}
(The key point here is that $V_{\attr}$ is is valid for the attractor \attr, not the tiling as in the previous theorem.) The precise \conds under which this may be accomplished are given in \cite{GeometryOfSST}. The forthcoming survey paper \cite{FractalTubeSurvey} gives a detailed discussion of the role of the present construction in applications related to tube formulas.
\section{The self-affine tiling} \label{sec:The Self-Similar
Tiling}

\subsection{Basic terms}
\label{ssec:Basic terms}

\begin{defn}
  \label{def:self-affine_system}
  A \emph{self-affine system} is a
  family $\simt := \{\simt_\j\}_{\j=1}^J$ (with $2 \leq J < \iy$) of
  affine mappings whose eigenvalues \gl all \sat $0< \gl< 1$. A
  \emph{self-similar system} is the special case where each mapping
  is a similitude, i.e.,
  \[\simt_\j(x) := r_\j A_\j x + a_\j,\]
  where for $\j=1,\dots,J$, we have $0 < r_\j < 1$, $a_\j \in \bRd$, and $A_\j
  \in O(d)$, the orthogonal group of $d$-\dimnl Euclidean space \bRd.
  The numbers $r_\j$ are referred to as the \emph{scaling ratios} of
  $\simt$. Thus, a similarity is a composition of an
  (affine) isometry and a homothety (scaling). All the examples
  presented here are self-similar, with the exception of the
  harmonic gasket, Example \ref{exm:harmonic_sg}.
\end{defn}

\begin{remark}
  Different self-affine systems may give rise to the same
  self-affine set. In this paper, the emphasis is placed on the
  self-affine system and its corresponding dynamical system, rather
  than on the self-affine set.
\end{remark}

\begin{defn}
  \label{def:attractor}
  A self-affine system is thus just a particular type of iterated
  function system (IFS). It is well known\footnote{See \cite{Hut}, as
  described in \cite{Fal1} or \cite{Kig}, for example.} that for such a
  family of maps, there is a unique and self-affine set \attr
  satisfying the fixed-point equation
  \begin{equation}
    \label{def:F}
    \attr = \simt(\attr) := \bigcup_{j=1}^J \simt_j(\attr).
  \end{equation}
  We call \attr the \emph{attractor} of \simt, or the
  \emph{self-affine set} associated with \simt. The \emph{action of}
  \simt is the set map defined by \eqref{def:F}. Thus, one says
  that \attr is invariant under the action of \simt.
\end{defn}

\begin{defn}
  \label{def:hull}
  We fix some notation for later use. Let
  \begin{equation}
    \label{eqn:hull def}
    \hull := [\attr]
  \end{equation}
  be the convex hull of the attractor \attr, that is, the collection of
  all convex combinations of points in \attr. (Equivalently, $[\attr]$ is
  the intersection of all convex sets containing \attr.) Since \attr is a
  \cpt set, it follows that \hull is also \cpt, by \cite[Thm.~1.1.10]{Schn2}.
  Further, let
  \begin{equation}
    \label{def:tileset}
    \hullint := \inter(\hull) = \hull \less \del \hull.
  \end{equation}
  Here, $\bdy A := \cj{A} \cap \cj{A^c}$, where $A^c$ is the complement
  of $A$ and $\cj A$ denotes the (topological) closure of $A$.
\end{defn}

\begin{remark}
  \label{rem:ambient_dim}
  For this paper, it will suffice to work with the ambient \dimn
  \begin{equation}
    \label{eqn:dim_of_hull}
    d = \dim \hull,
  \end{equation}
  restricting the maps $\simt_j$ as appropriate. In
  \eqref{eqn:dim_of_hull}, $\dim \hull$ is defined to be the usual
  topological \dimn of the smallest affine space containing \hull.
  An appropriate change of coordinates allows one to think of this
  convention as using a minimal space \bRd in which to embed \attr;
  e.g., if \attr is a Cantor set in \bRT, we study it as if the ambient
  space were the line containing it, rather than \bRT. Note that this means
  \hullint is open in the standard topology; and so we have $\hull^o
  \neq \es$. This remark is intended to allay any fears about
  possibly needing to use relative interior instead of interior (see
  \cite{Rota} or \cite{Schn2}) and other unnecessary complications.
\end{remark}

\begin{defn}
  \label{def:tileset_condition}
  A self-affine system \sats the \emph{tileset condition} iff
  for $j \neq \ell$,
  \begin{equation}
    \label{eqn:tileset_condition}
    \inter \simt_j(\hull) \cap \inter \simt_\ell(\hull) = \es.
  \end{equation}
  It is shown in Cor.~\ref{cor:tileset_implies_boundarylaps}
  that \bc $\hull = \cj{\inter \hull}$,
  \eqref{eqn:tileset_condition} implies that the images
  $\simt_j(\hull)$ and $\simt_\ell(\hull)$ can intersect only on
  their boundaries:
  \begin{equation*}
    \simt_j(\hull) \cap \simt_\ell(\hull)
    = \bdy \simt_j(\hull) \cap \bdy \simt_\ell(\hull).
  \end{equation*}
  To avoid trivialities, we also require
  \begin{equation}
    \label{eqn:tileset_nontriviality}
    \hullint \nsubseteq \simt(\hull).
  \end{equation}
  The \emph{nontriviality condition} \eqref{eqn:tileset_nontriviality}
  disallows the case $\hullint \less \simt(\hull)= \es$, and hence
  guarantees the existence of the tiles in \S\ref{ssec:The_construction}. It is shown in \cite{GeometryOfSST} that \eqref{eqn:tileset_nontriviality} fails iff \attr is convex.
\end{defn}

\begin{remark}
  \label{rem:tileset_implies_OSC}
  The tileset condition is a restriction on the overlap of the images
  of the mappings, comparable to the \emph{open set condition} (OSC).
  The OSC requires a nonempty bounded open set $U$ \st the sets $\simt_j(U)$
  are disjoint but $\simt(U) \ci U$. See, e.g., \cite[Chap.~9]{Fal1}.
  It is clear from Cor.~\ref{cor:hull_k traps A} of \S\ref{sec:Properties of the tiling}
  (with $k=0$ and $U=\inter\hull$) that the OSC follows from
  \eqref{eqn:tileset_condition}. However, the following example
  shows that the converse is false.
\end{remark}

\begin{figure}
  \centering
  \scalebox{0.85}{\includegraphics{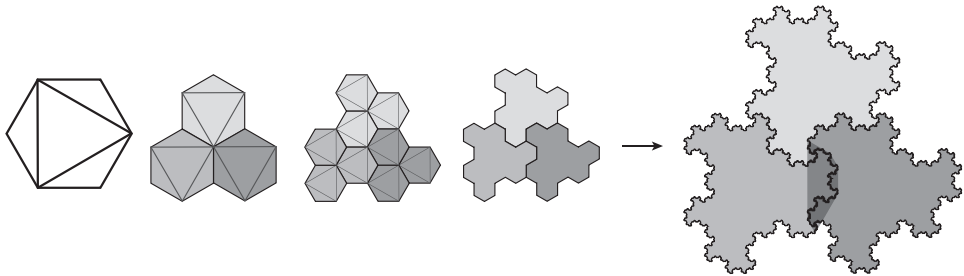}}
  \caption{\captionsize A self-similar system which satisfies the
  open set condition but not the tileset condition; see
  Example~\ref{exm:tileset_counterexample}. The attractor in
  this example tiles all of \bRt.
  }
  \label{fig:tileset_counterexample}
\end{figure}

\begin{ex}
  \label{exm:tileset_counterexample}
  Consider a system of three similarity mappings, each with scaling ratio
  $1/\sqrt3$ and a clockwise rotation of $\gp/2$. The mappings
  are illustrated in Figure~\ref{fig:tileset_counterexample} and
  form a system which \sats the open set condition (simply take
  the interior of the attractor) but not the tileset condition.
  On the right, the attractor
  has been shaded for clarity; the dark overlay indicates the
  intersection of the convex hulls of the lower two images of the
  attractor.
\end{ex}

\begin{defn}
  Denote the \emph{words of length $k$} (of
  $\{1,2,\dots,J\}$) by
  \begin{align}
    \label{eqn:Wk_of_J}
    W_k = W_k^J :&= \{1,2,\dots,J\}^k \notag \\
           &= \{w=w_1 w_2 \dots w_k \suth w_j \in \{1,2,\dots,J\}\},
  \end{align}
  and the set of all (finite) words by
  \(
    W := \bigcup_k W_k\).
  Generally, the dependence of $W_k^J$ on $J$ is suppressed.
  For $w$ as in \eqref{eqn:Wk_of_J}, we use the standard IFS notation
  \begin{equation}
    \label{def:comp_via_w}
    \simt_w(x) := \simt_{w_1} \comp \simt_{w_2} \comp \dots \comp \simt_{w_k}(x)
  \end{equation}
  to describe compositions of maps from the self-affine system.
\end{defn}

\begin{defn}
  \label{def:tiling}
  Let $A \ci \bRd$ be a set which is the closure of its interior.
  A \emph{tiling} of $A$ is a collection of nonempty connected $d$-\dimnl
  sets $\{A_n\}_{n=1}^\iy$ \st
  \begin{Cond}
    \item $A = \cj{\bigcup_{n=1}^\iy {A_n}}$, and
    \item $A_n \cap A_m \ci \del A_n \cap \del A_m$ for $n \neq m$.
  \end{Cond}
  We then say that the sets $A_n$ \emph{tile} $A$. Further, define
\end{defn}

\begin{figure}
  \scalebox{1.0}{\includegraphics{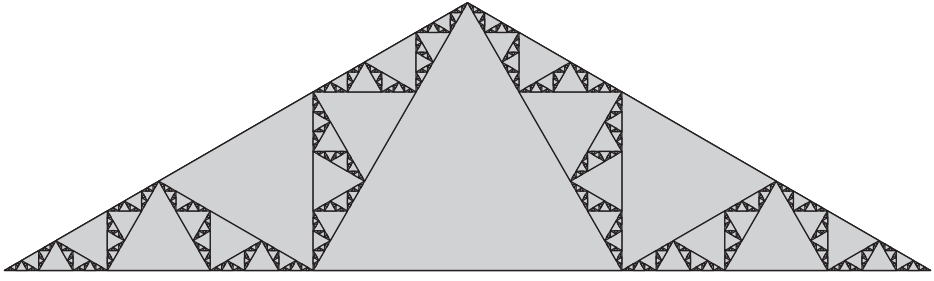}}
  \caption{\captionsize Tiling the complement of the Koch curve $K$. The equilateral
  triangles form an open tiling of the convex hull $\hull = [K]$, in
  the sense of Def.~\ref{def:tiling}.
  }
  \label{fig:koch-tiled-in}
  \centering
\end{figure}

\begin{defn} \label{def:open_tiling}
  A \emph{tiling of $A$ by open sets} (or \emph{open tiling})
  is a collection of nonempty connected open sets $\{A_n\}_{n=1}^\iy$ \st
  \begin{Cond}
    \item $\cj{A} = \cj{\bigcup_{n=1}^\iy {A_n}}$, and
    \item $A_n \cap A_m = \es$ for $n \neq m$.
  \end{Cond}
  Figure~\ref{fig:koch-tiled-in} shows an example of a tiling by open
  sets, each of which is an equilateral triangle.
\end{defn}

\subsection{The construction}
\label{ssec:The_construction}

In this section, we construct a \emph{self-affine tiling}, that is,
a tiling which is constructed via the action of a self-affine system
(as defined in Def.~\ref{def:self-affine_system}). Such a tiling
will consequently have a self-affine structure, and is defined
precisely in Def.~\ref{def:ss tiling} below. The reader is invited
to look ahead at Figure~\ref{fig:koch-contractions}, where the
construction is illustrated step-by-step for the illuminative
example of the Koch curve.

For the system $\{\simt_j\}$ with attractor \attr, denote the convex
hull of the attractor by
\begin{equation}
  \label{eqn:hull0}
  \hull_0 = \hull := [\attr].
\end{equation}
Denote the image of \hull under the action of \simt (in accordance
with \eqref{def:F}) by
\begin{equation}
  \label{eqn:hullk}
  \hull_k := \simt^k(\hull) = \bigcup_{w \in W_k} \simt_w(\hull), \q k=1,2,\dots \,.
\end{equation}
Note that this is equivalent to the inductive definition
\begin{equation}
  \label{eqn:hullk_inductively}
  \hull_k := \simt(\hull_{k-1}), \q k=1,2,\dots \,.
\end{equation}

\begin{defn}
  \label{def:tilesets}
  The \emph{tilesets} are the sets
  \begin{align}
    \tileset_k &:= \cj{\hull_{k-1} \less \hull_k}, \q k=1,2,\dots \,
      \label{eqn:tilesetk}
  \end{align}
\end{defn}

\begin{defn}
  \label{def:generators}
  The \emph{generators} $\gen_q$ are the connected components of the open set
  \begin{equation}
    \label{eqn:tileset1 as union}
    \inter (\hull \less \simt(\hull))
      =\gen_1 \sqcup \gen_2 \sqcup \dots \sqcup \gen_Q.
  \end{equation}
\end{defn}

\begin{remark}
  \label{rem:disjoint_union}
  The symbol $\sqcup$ is used in place of $\cup$ to emphasize that a
  given union of sets is \emph{disjoint}. This should not be
  confused with the operation of \emph{disjoint union}, i.e., the
  coproduct in the category of sets.
\end{remark}

\begin{figure}
  \centering
  \scalebox{1.05}{\includegraphics{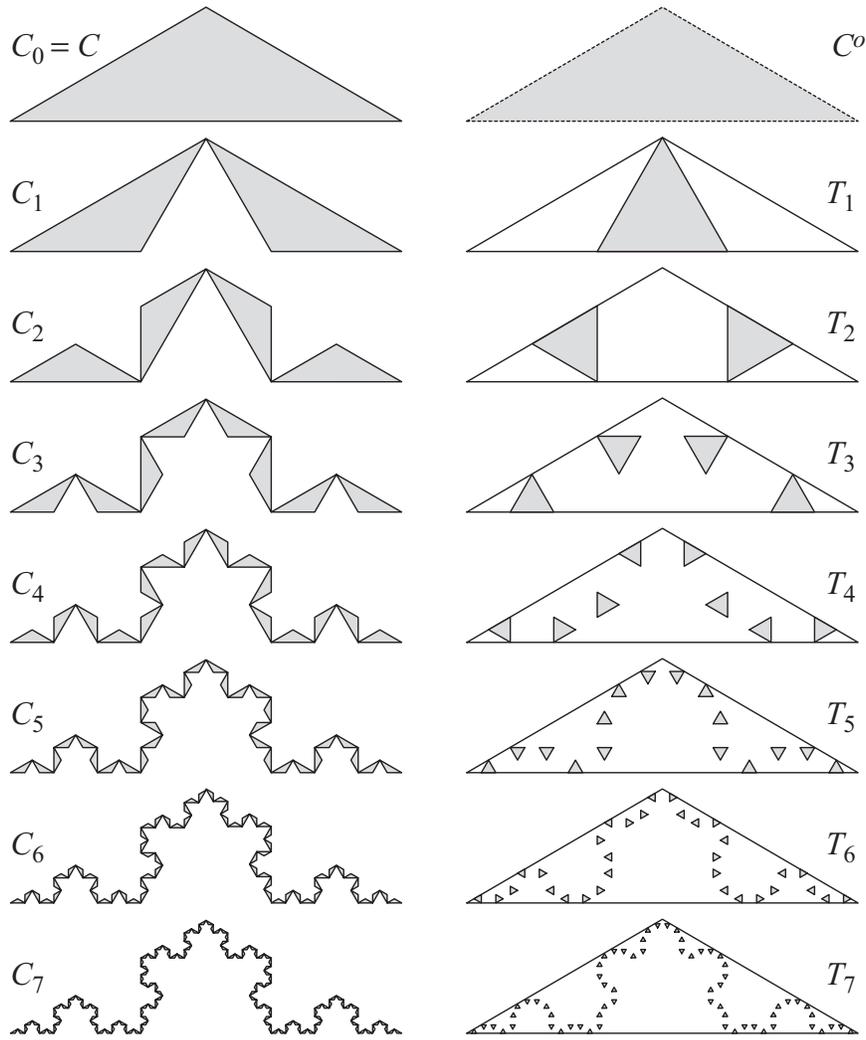}}
  \caption[Tiling the Koch curve; contractions of the convex hull.]
    {The left column shows images of the convex hull
  $\hull$ under successive applications of \simt. The right
  column shows how the components of the $\tileset_k$ tile the
  complement; they are overlaid in Figure~\ref{fig:koch-tiled-in}.
  This tiling has one generator \mbox{$\gen_1= \inter\tileset_1$}.}
  \label{fig:koch-contractions}
\end{figure}

As will be shown in Theorem~\ref{thm:tilesets_are_nondegenerate}, it
follows from the tileset conditions
\eqref{eqn:tileset_condition}--\eqref{eqn:tileset_nontriviality}
(and some other facts) that the tilesets and tiles are nonempty, and
that each tileset is the closure of its interior. Also,
Theorem~\ref{thm:propagation_of_tiles} will justify the terminology
``generators'' by showing
\begin{equation}
  \label{eqn:tileset1 as union}
  \tileset_{k} = \cj{\bigsqcup\nolimits_{q=1}^Q \simt^{k-1}(\gen_q)},
\end{equation}
that is, that any difference $\hull_{k-1} \less \hull_k$ is (modulo
some boundary points) the image of the generators under the action
of \simt. The number $Q$ of generators depends on the specific
geometry of $\hull$ and on the self-affine system \simt. It is
conceivable that $Q=\iy$ for some systems \simt, but no such
examples are known. This possibility will be investigated further in
\cite{GeometryOfSST}.

\begin{defn}
  \label{def:ss tiling}
  The \emph{self-affine tiling} of \attr is
  \begin{equation}
  \label{eqn:ss tiling}
    \sT := \left(\{\simt_j\}_{j=1}^J,\{\gen_q\}_{q=1}^Q\right).
  \end{equation}
\end{defn}

We may also abuse the notation a little, and use \tiling to denote
the set of corresponding \emph{tiles}:
\begin{equation}
  \label{def:generators}
  \tiling = \{\tile_n\}_{n=1}^\iy
  = \{\simt_w(\gen_q) \suth w \in W, q=1,\dots,Q\},
\end{equation}
where the \seq $\{\tile_n\}$ is an enumeration of the tiles.
Clearly, each tile is nonempty and $d$-\dimnl. Furthermore,
Theorem~\ref{thm:the tiling} will confirm that
\eqref{def:generators} is a tiling by open sets, as in
Def.~\ref{def:tiling}.

\section{Examples} \label{ssec:Examples}

All the examples discussed in this section have polyhedral
generators, but this is not the general case. In fact, it is
possible to have generators with boundary that is continuously
differentiable, although it is not possible that they be twice
continuously differentiable. This was observed to be true for the
convex hull of an attractor in \cite{GAST1}, and it immediately
carries over to the generators as well. We will study this
eventuality further in \cite{GeometryOfSST} and \cite{TFG}. See also
Remark~\ref{ssec:Properties of the generators}.

\subsection{The Koch curve.}
\label{exm:koch}

Figure~\ref{fig:koch-tiled-in} shows the self-affine tiling of the
Koch curve; the steps of the construction are illustrated in
Figure~\ref{fig:koch-contractions}. The tiling is \(\sK =
\left(\{\simt_j\}_{j=1}^2,\{\gen\}\right)\), and it is easiest to
write down the maps as $\simt_j:\bC \to \bC$, with the natural
identification of \bC and \bRt. In this case, $\simt_1(z) := \gx
\cj{z}$ and $\simt_2(z) := (1-\gx)(\cj{z}-1)+1.$
Figure~\ref{fig:koch-contractions} depicts the case when $\gx = (1 +
\sqrt{-1/3})/2$, so that $r_1 = r_2 = 1/\sqrt3$ and the single
generator $\gen$ is the equilateral triangle of side length
$\tfrac13$.

In general, \gx may be any complex number \satg $|\gx|^2 + |1-\gx|^2
< 1,$ i.e., lying within the circle of radius 1/2 centered at
$z=1/2$. Basic geometric considerations show that this inequality
must be satisfied in order for the tileset condition
\eqref{eqn:tileset_condition} to be met. Any member of this family
will have one isosceles triangle $\gen=\gen_1=\tileset_1$ for a
generator. A key point of interest in this example is that, in the
language of \cite{FGCD}, curves from this family will generically be
nonlattice, i.e., the logarithms of the scaling ratios will not be
rationally dependent. Figure~\ref{fig:koch-contractions} thus shows
a very exceptional case.

\subsection{The Sierpinski gasket} \label{exm:Sierpinski gasket}

\begin{figure}[b]
    \centering
    \scalebox{1.0}{\includegraphics{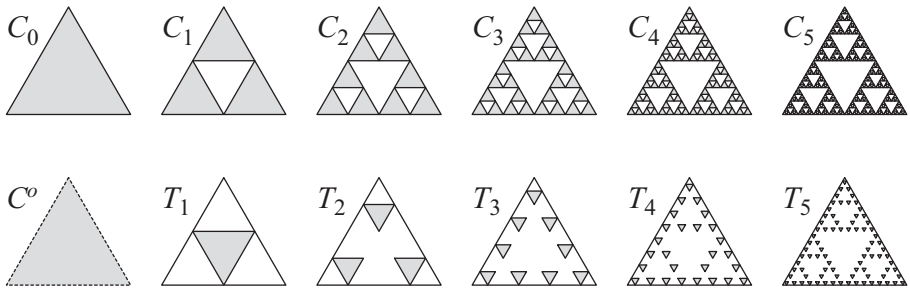}}
    \caption[Self-similar tiling of the Sierpinski gasket.]
      {\captionsize Self-similar tiling of the Sierpinski gasket.}
    \label{fig:sierpinski-gasket-contractions}
  \end{figure}

  The Sierpinski gasket system consists of the three maps
  $\simt_j(x) = (x + p_j)/2,$
  where the $p_j$ are the vertices of an equilateral triangle;
  the standard example is
  \(p_1 = (0,0), p_2 = (1,0), \text{ and } p_3 = (\tfrac12,\frac{\sqrt3}2.\)
  The convex hull of the gasket is the triangle with
  vertices $p_1,p_2,p_3$. The generator \gen is the `middle
  fourth' of the hull (see $\tileset_1$ in
  Figure~\ref{fig:sierpinski-gasket-contractions}).

  \begin{figure}[b]
    \centering
    \scalebox{0.80}{\includegraphics{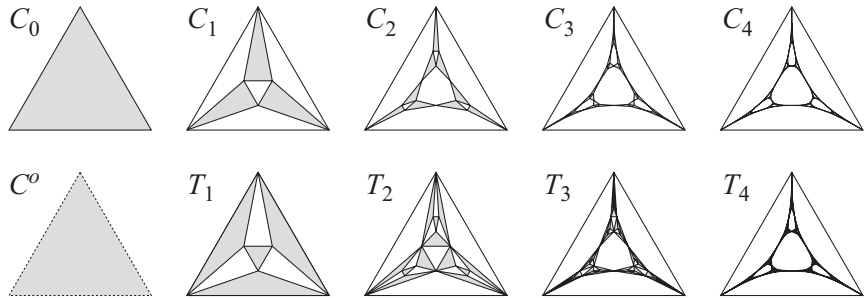}}
    \caption[Self-affine tiling of the harmonic Sierpinski gasket.]
      {\captionsize Self-affine tiling of the harmonic Sierpinski gasket.}
    \label{fig:harmonic-gasket-contractions}
  \end{figure}

\subsection{The harmonic Sierpinski gasket.}
\label{exm:harmonic_sg}

Recent studies of PDE on post-critically finite fractals (following
the analytic methods of \cite{Kig}) have led to interest in the
Sierpinski gasket as it is embedded into \bRt via harmonic
coordinates. The resulting figure is a self-affine homeomorphic
image of the usual gasket which is not self-similar. The eigenvalues
of the affine maps are $\frac35,\frac15$. See
Figure~\ref{fig:harmonic-gasket-contractions}.

\subsection{The pentagasket} \label{exm:pentagasket}

\begin{figure}
    \centering
    \scalebox{1.15}{\includegraphics{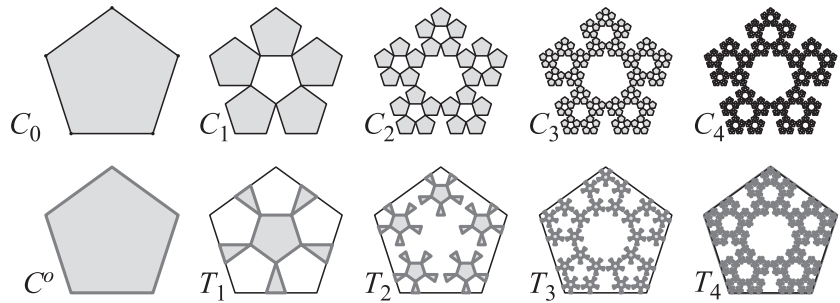}}
    \caption[Self-affine tiling of the pentagasket.]
      {\captionsize Self-affine tiling of the pentagasket.}
    \label{fig:pentagasket}
  \end{figure}

  The pentagasket is constructed via five maps
  $\simt_j(x) = \phi^{-2} x + p_j,$
  where the $p_j$ form the vertices of a pentagon of side length 1,
  and $\phi = \frac{1+\sqrt5}2$ is the golden ratio, so that the
  scaling ratio of each mapping is
  \[r_j = \phi^{-2} = \tfrac{3-\sqrt5}2, \q j=1,\dots,5.\]
  The pentagasket is a self-similar (not just self-affine) example of multiple generators
  $\gen_q$ with $q=1,\dots,6$. In fact, $\tileset_1=\gen_1 \cup \dots \cup \gen_6$
  where $\gen_1$ is a pentagon and $\gen_2,\dots,\gen_6$ are
  triangles. See Figure~\ref{fig:pentagasket}.

\subsection{The Sierpinski carpet} \label{exm:Sierpinski carpet}

\begin{figure}
    \centering
    \scalebox{1.0}{\includegraphics{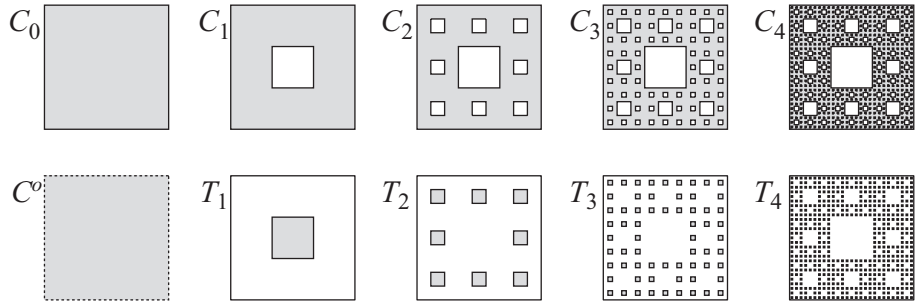}}
    \caption[Self-similar tiling of the Sierpinski carpet.]
      {\captionsize Self-similar tiling of the Sierpinski carpet.}
    \label{fig:sierpinski-carpet-contractions}
  \end{figure}

  The Sierpinski carpet is constructed via eight maps
  $\simt_j(x) = (x + p_j)/3,$
  where $p_j = (a_j,b_j)$ for $a_j,b_j\in \{0,1,2\}$,
  excluding the single case $(1,1)$. The Sierpinski carpet is an
  example which is not finitely ramified; indeed, it is not even
  post-critically finite (see \cite{Kig}).  See Figure~\ref{fig:sierpinski-carpet-contractions}.

\subsection{The Menger sponge} \label{exm:Menger sponge}
\begin{figure}
    \centering
    \scalebox{1.0}{\includegraphics{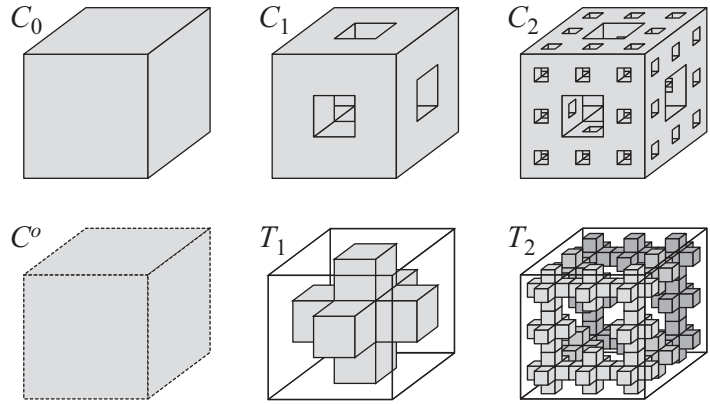}}
    \caption[Self-similar tiling of the Menger sponge.]
      {\captionsize Self-similar tiling of the Menger sponge.}
    \label{fig:menger-sponge-contractions}
  \end{figure}

  The Menger sponge is constructed via twenty maps
  $\simt_j(x) = (x + p_j)/3,$ where $p_j = (a_j,b_j,c_j)$ for
  $a_j,b_j,c_j \in \{0,1,2\}$,
  except for the six cases when exactly two coordinate are $1$, and the
  single case when all three coordinates are $1$.
  The Menger sponge system provides an example with an generator of
  dimension greater than 2, and is also an example with a nonconvex
  generator. See Figure~\ref{fig:menger-sponge-contractions}.

\pgap

\section{Properties of the tiling}
\label{sec:Properties of the tiling}

The results of this section indicate that a self-affine tiling may
be constructed for any self-affine system \satg the tileset
condition (and nontriviality condition) of
Def.~\ref{def:tileset_condition}. Throughout, we will
use the fact that $\cj{A} = \inter A \sqcup \bdy A$, 
where we denote the closure of $A$ by $\cj A$, the interior of $A$
by $\inter A$, and the boundary of $A$ by $\bdy A = \cj{A} \cap
\cj{A^c}$, where $A^c = \bRd \less A$. Recall from
Rem.~\ref{rem:disjoint_union} that $\sqcup$ indicates a union of
\emph{disjoint} sets.

The first part of this section,
Thm.~\ref{thm:successor_is_always_contained}--Cor.~\ref{cor:hulls_converge_to_F},
establishes the nested structure of the attractor and the images of
the hull \hull,
\[\hull_{k+1} \ci \hull_k, \text{ and } \bigcap \hull_k = \attr.\]
These results are reminiscent of \cite[5.2(3)]{Hut}, but are
developed in terms of the convex hull of the attractor, rather than
a set satisfying the open set \cond.

\begin{theorem}
  \label{thm:successor_is_always_contained}
  For each $k \in \bN$, one has $\hull_{k+1} \ci \hull_k \ci \hull$.
  \begin{proof}
    Any point $x \in \hull$ is a convex combination of points in
    \attr. Since affine transformations preserve convexity,
    $\simt_j(x)$ will be a convex combination of points in
    $\simt_j(\attr) \ci \attr$. Hence $\simt_j(\hull) \ci [\attr]
    = \hull$ for each $j$, so $\simt(\hull) \ci \hull$. By iteration
    of this argument, we immediately have $\simt^k(\hull) \ci \hull$
    for any $k \in \bN$.
    From \eqref{eqn:hullk_inductively}, it is clear that
    \begin{align}
      \label{eqn:successor_containment}
      \hull_{k+1}
        = \simt(\hull_{k})
        = \simt^{k+1}(\hull)
        = \simt^k(\simt(\hull))
        \ci \simt^k(\hull)
        = \hull_k,
    \end{align}
    where the inclusion follows by $\simt(\hull) \ci \hull$, as
    established initially.
  \end{proof}
\end{theorem}

\begin{cor}
  \label{cor:tileset_condition_iterates}
  The tileset condition is preserved under the action of \simt, i.e.,
  \begin{equation}
    \label{eqn:tileset_condition_iterates}
    \inter \simt_j(\hull_k) \cap \inter \simt_\ell(\hull_k) = \es,
    \q j \neq \ell,
    \q \forall k \in \bN.
  \end{equation}
  \begin{proof}
    From Theorem~\ref{thm:successor_is_always_contained} we have
     $\inter \simt_j(\hull_k) \ci \inter \simt_j(\hull)$,
    and similarly for $\simt_\ell$.
    Then \eqref{eqn:tileset_condition_iterates} follows immediately
    from the tileset \cond \eqref{eqn:tileset_condition}.
  \end{proof}
\end{cor}

\begin{cor}
  \label{cor:hull_k traps A}
  For $A \ci \hull_k$, we have $\simt_w(A) \ci \hull_k$, for all $w \in W$.
  In particular, $\attr \ci \hull_k, \forall k$.
  \begin{proof}
    By iteration of \eqref{eqn:successor_containment},
    it is immediate that $\hull_m \ci \hull_k$ for any $m \geq k$.
    Since $\simt(A) \ci \simt(\hull_k) = \hull_{k+1} \ci \hull_k$ by
    Theorem~\ref{thm:successor_is_always_contained}, the first
    conclusion follows.
    The special case follows by
    induction on $k$ with basis case $A = \attr \ci \hull = \hull_0$.
    The inductive step is
    \[\attr \ci \hull_k
      \q \implies \q
      \attr = \simt(\attr)
        \ci \simt(\hull_k)
        = \hull_{k+1}.
        \qedhere\]
  \end{proof}
\end{cor}

\begin{remark}
  Cor.~\ref{cor:hull_k traps A} is a modified form of
  \cite[3.1(8)]{Hut}.
\end{remark}

\begin{cor}
  \label{cor:hulls_converge_to_F}
  The decreasing \seq of sets $\{\hull_k\}$ converges to \attr.
  \begin{proof}
    Cor.~\ref{cor:hull_k traps A} shows
    $\attr \ci \hull_k$ for every $k$, so it is clear that
    $\attr \ci \bigcap \hull_k$.
    For the reverse inclusion, suppose $x \notin \attr$, so that
    $x$ must be some positive distance \ge from \attr. Let
    $\gl$ be the largest eigenvalue of the maps $\{\simt_j\}$,
    and recall that $0 < \gl < 1$. For $w \in \Wds_k$, we have
    $\diam(\simt_w(\hull)) \leq \gl^k\diam(\hull)$, which clearly
    tends to 0 as $k \to \iy$.
    Therefore, we can find $k$ for which all points of $\hull_k =
    \simt^k(\hull)$ lie within $\ge / 2$ of \attr.
    Thus $x$ cannot lie in $\hull_k$ and hence $x \notin \bigcap
    \hull_k$.
  \end{proof}
\end{cor}

\begin{remark}
  \label{rem:hulls_converge_to_F_in_Hausdorff}
  The convergence $\hull_k \to \attr$ also holds in the sense of Hausdorff
  metric, by \cite[3.2(1)]{Hut}; see also \cite{Fal1} or \cite{Kig} for a nice
  discussion. Hutchinson showed that \simt is a contraction mapping on
  the metric space of compact subsets of \bRd, which is complete when
  endowed with the Hausdorff metric. (The hull \hull is shown to be
  compact in Cor. \ref{cor:Cks_are_nondegenerate}.)
  An application of the contraction
  mapping principle then shows that \simt has a unique attracting
  fixed point (as stated in Def.~\ref{def:attractor}). This phenomenon
  is apparent in several of the figures. The main observation in the
  current situation is that $\hull_k$ \emph{decreases} to \attr, by
  the nestedness described in
  Theorem~\ref{thm:successor_is_always_contained}.
\end{remark}

We now use the conditions on the mapping system \simt to obtain some
technical nondegeneracy results in
Lemma~\ref{lemma:simt_preserves_closures}--Cor.~\ref{cor:simtj_preserves_closint}.
There are two main ideas:
\begin{enumerate}
  \item All the sets we work with ($\hull_k, \tileset_k$, etc.) are
    nondegenerate in the sense of being the closure of their interior.
  \item The action of \simt preserves several key properties,
    including closure, the tileset conditions, and the nondegeneracy
    condition mentioned just above.
\end{enumerate}
In the latter part of this section,
Thm.~\ref{thm:propagation_of_tilesets}--Cor.~\ref{cor:tiling_is_subselfaffine},
these allow us to connect properties of the hull differences
$\hull_k \less \hull_{k+1}$ to properties of the tilesets
$\tileset_k$. Heuristically, we are leveraging the relation
\[\hull_k \less \hull_{k+1} \;\approx\; \simt^k(\gen_1 \cup \dots \cup \gen_q)\]
to show that the sets $\simt_w(\gen_q)$ give a tiling of \hull by
open sets, i.e., the construction can always be carried out for an
iterated \fn system satisfying the tileset condition. The displayed
equation above is only approximate because the sets may differ at
the boundary; unfortunately, this necessitates several
technicalities in the proofs of the main results.

\begin{lemma}
  \label{lemma:simt_preserves_closures}
  The action of \simt commutes with set closure, i.e., $\simt\left(\cj A\right) = \cj{\simt(A)}$
  \begin{proof}
    It is well known that closure commutes with finite unions, i.e.,
    for any sets $A,B$, one has $\cj A \cup \cj B = \cj{A \cup B}$.
    See, e.g., \cite[Chap.~2, \S17]{Munk}. Also, each $\simt_j$ is a
    \homeo, and is thus a closed, \cn map. Therefore,
    \[\simt\left(\cj A\right)
      = \bigcup\nolimits_{j=1}^J \simt_j\left(\cj A\right)
      = \bigcup\nolimits_{j=1}^J \cj{\simt_j(A)}
      = \cj{ \bigcup\nolimits_{j=1}^J \simt_j(A)}
      = \cj{ \simt(A)}.
      \qedhere\]
  \end{proof}
\end{lemma}


\begin{theorem}
  \label{thm:simt_preserves_nondegeneracy}
  If $A$ is the closure of its interior, then so is $\simt(A)$.
  \begin{proof}
    Since $\simt_j$ is a \homeo, this is a simple exercise in basic topology.
  \end{proof}
\end{theorem}


\begin{cor}
  \label{cor:Cks_are_nondegenerate}
  Each set $\hull_k$ is the closure of its interior.
  \begin{proof}
    The set $\hull = [\attr]$ is convex by definition, and compact
    by \cite[Thm.~1.1.10]{Schn2}. Therefore, \hull is the closure
    of its interior by \cite[Thm.~1.1.14]{Schn2}. The conclusion
    follows by iteration of Theorem~\ref{thm:simt_preserves_nondegeneracy}.
  \end{proof}
\end{cor}


\begin{cor}
  \label{cor:tileset_implies_boundarylaps}
  The tileset condition implies that images of the hull can only
  overlap on their boundaries:
  \begin{align}
    \simt_j(\hull) \cap \simt_\ell(\hull)
    = \bdy \simt_j(\hull) \cap \bdy \simt_\ell(\hull),
    \q \text{ for } j \neq \ell.
  \end{align}
  \begin{proof}
    Let $x \in \simt_j(\hull) \cap \simt_\ell(\hull)$. Suppose,
    by way of contradiction, that $x \in \inter \simt_j(\hull)$.
    Then we can find an open \nbd $U$ of $x$ which is contained in
    $\inter \simt_j(\hull)$. Since $x \in \simt_\ell(\hull)$,
    there must be some $z \in U \cap \inter \simt_\ell(\hull)$, by
    Cor.~\ref{cor:Cks_are_nondegenerate}. But then $z \in \inter
    \simt_j(\hull) \cap \inter \simt_\ell(\hull)$, in contradiction
    to the tileset condition. For the reverse inclusion, note that
    $\del A \cap \del B \ci A \cap B$ whenever $A, B$ are closed
    sets.
  \end{proof}
\end{cor}

\begin{theorem}[Nondegeneracy of tilesets]
  \label{thm:tilesets_are_nondegenerate}
  Each tileset is the closure of its interior.
  \begin{proof}
    We need only show $\tileset_k \ci \cj{\inter \tileset_{k}}$,
    since the reverse containment is clear by the closedness of
    $\tileset_k$. Since $\cj A = \inter A \sqcup \bdy A$,
    take $x \in \inter (\hull_{k-1} \less \hull_{k})$ to begin.
    Using Cor.~\ref{cor:Cks_are_nondegenerate}, we have equality
    in the first step of the following derivation:
    \begin{align}
      \hull_{k-1} \less \hull_k
      = \cj{\inter(\hull_{k-1})} \less \cj{\hull_k}
      &\ci \cj{\inter(\hull_{k-1}) \less \hull_k} 
      \ci \cj{\inter\left(\cj{\hull_{k-1} \less \hull_k} \right)},
    \end{align}
    \bc 
    $\inter(\hull_{k-1}) \less \hull_k
      = \inter(\inter(\hull_{k-1}) \less \hull_k)
      \ci \inter(\hull_{k-1} \less \hull_k).
    $

    Now consider the case when $x \in \bdy(\hull_{k-1} \less \hull_k)$.
    Pick an open \nbd $U$ of $x$. By the same argument as above,
    choose $z \in U \cap (\hull_{k-1} \less \hull_k)$ to see that
    $x$ is a limit point (and hence an element) of $\cj{ \inter \tileset_k}$.
  \end{proof}
\end{theorem}

\pgap

The following two corollaries will be useful in the proof of
Theorem~\ref{thm:propagation_of_tiles}.

\begin{cor}
  \label{cor:simtj_preserves_closint}
  For $j=1,\dots,J$, $\cj{\simt_j(\hull_{k-1}) \less \simt_j(\hull_{k})}$
  is the closure of its interior.
  \begin{proof}
    Because each $\simt_j$ is a \homeo, the set $\simt_j(\cj{\hull_{k-1}
    \less \hull_{k}})$ will be the closure of its interior by
    Theorem~\ref{thm:tilesets_are_nondegenerate}. However, we have
    \begin{align}
      \simt_j(\cj{\hull_{k-1} \less \hull_{k}})
      = \cj{\simt_j(\hull_{k-1} \less \hull_{k})}
      = \cj{\simt_j(\hull_{k-1}) \less \simt_j(\hull_{k})},
        \label{eqn:simtj_preserves_diffs_and_clos}
    \end{align}
    since $\simt_j$ is closed and injective.
  \end{proof}
\end{cor}

\begin{cor}
  \label{cor:hull-less-simt(hull)-is-nice}
  It is always the case that one has $\cj{\hull \less \simt(\hull)} = \cj{\inter(\hull \less \simt(\hull))}$
  and $\del(\hull \less \simt(\hull)) = \del(\inter(\hull \less \simt(\hull)))$.
\end{cor}
\begin{proof}
  For the first statement, $\inter(\hull \less \simt(\hull)) \ci
  \hull \less \simt(\hull)$ immediately shows one
  containment. For the other, choose $x \in \cj{\hull \less
  \simt(\hull)}$ so that there is a \seq $x_k \to x$ with $x_k
  \in \hull \less \simt(\hull)$. Since $\simt(\hull) =
  \cj{\simt(\hull)}$, we have $x_k \in \cj{\simt(\hull)}^c =
  \inter(\simt(\hull)^c)$. Further, since $x_j \in \hull =
  \cj{\inter(\hull)}$, one can always find
  \[y_k \in B(x_k,\tfrac1k) \cap \inter(\hull) \cap
    \inter(\simt(\hull)^c).\]
  Then we have a \seq $y_k \to x$, and since
  $\inter(\simt(\hull)^c) = \cj{\simt(\hull)}^c =
  \simt(\hull)^c$, one has
  \begin{align*}
    y_k \in \inter(\hull) \cap \inter(\simt(\hull)^c)
    &= (\inter \hull) \less \simt(\hull) \\
    &= \inter(\inter \hull \less \simt(\hull)) \\
    &\ci \inter(\hull \less \simt(\hull)).
  \end{align*}
  Thus, $x \in \cj{\inter(\hull \less \simt(\hull))}$. To see
  the second statement holds, it suffices to show $\cj{(\hull \less \simt(\hull))^c}
  = \cj{(\inter(\hull \less \simt(\hull)))^c}$. However, this
  is clear by applying the identity $(\inter A)^c = \cj{A^c}$
  to the identity $\inter(\hull \less \simt(\hull)) =
  \inter(\inter(\hull \less \simt(\hull)))$.
\end{proof}

\pgap

We are now ready to prove the main results of this paper. The reader
may find Figure \ref{fig:koch-tiled-in} and Figure
\ref{fig:koch-contractions} helpful, as illustrations of
Theorems~\ref{thm:propagation_of_tilesets}--\ref{cor:tiling_is_subselfaffine}.
The next result shows that each tileset is the image under \simt of
its predecessor; speaking very roughly, the hulls $\hull_k$ form a
\seq of ``\nbds'' of the attractor \attr. In the sense of dynamical
systems, the tiles describe the orbits of \simt.

\begin{theorem}
  \label{thm:propagation_of_tilesets}
  Each tileset is the image under \simt of its predecessor, i.e.,
  \begin{align}
    \simt(\tileset_k) = \tileset_{k+1}, \q \text{for } k \in \bN.
      \label{eqn:propagation_of_tilesets}
  \end{align}
  \begin{proof}
    One would like to say simply that
    \begin{align*}
      \simt(\tileset_k) = \simt\left(\cj{\hull_{k-1} \less \hull_k}\right)
      = \cj{\simt\left(\hull_{k-1} \less \hull_k\right)}
      = \cj{\simt(\hull_{k-1}) \less \simt(\hull_k)}
      = \cj{\hull_{k} \less \hull_{k+1}}
      = \tileset_{k+1}.
    \end{align*}
    Unfortunately, the central equality is not immediately justifiable.
    Using using Def.~\ref{def:tilesets} and \eqref{eqn:simtj_preserves_diffs_and_clos},
    we have the identities
    \begin{align}
      \simt(\tileset_k)
      &= \bigcup\nolimits_{j=1}^J \cj{\simt_j(\hull_{k-1} \less \hull_k)}, \text{ and }
        \label{eqn:use_pres_of_inter} \\
      \tileset_{k+1}
      &= \cj{\hull_{k} \less \hull_{k+1}}.
        \label{eqn:use_of_tileset}
    \end{align}

    ($\ci$) To see that \eqref{eqn:use_pres_of_inter} is a subset of
    \eqref{eqn:use_of_tileset}, pick $x \in \simt(\tileset_k)$ and proceed
    by cases.

    (i) For any $x \in \inter(\simt_j(\hull_{k-1}) \less \simt_j(\hull_k))$,
    it must be that $x \in \inter \simt_j(\hull_{k-1}) \ci \hull_k$.

    It suffices to show $x \in \hull_{k} \less \hull_{k+1}$, so by way of contradiction,
    suppose that $x \in \hull_{k+1}$. Then $x \in
    \simt_\ell(\hull_k)$ for some $\ell \neq j$. Inasmuch as
    Theorem~\ref{thm:successor_is_always_contained} gives $x
    \in \simt_\ell(\hull_{k-1})$, Cor.~\ref{cor:tileset_implies_boundarylaps}
    implies $x \in \bdy \simt_j(\hull_{k-1}) \cap \bdy
    \simt_\ell(\hull_{k-1}),$
    contradicting the fact that $x \in \inter \simt_j(\hull_{k-1})$.

    (ii) Now consider $x \in \bdy (\simt_j(\hull_{k-1}) \less \simt_j(\hull_k))$.
    Let $U$ be an open \nbd of $x$. By
    Cor.~\ref{cor:simtj_preserves_closint}, we can find $w \in U \cap
    \inter(\simt_j(\hull_{k-1}) \less \simt_j(\hull_k))$. Repeating
    the arguments of part (i), we obtain $w \in \cj{\hull_k \less
    \hull_{k+1}}$. Thus $x$ is a limit point (and hence a member)
    of $\cj{\hull_k \less \hull_{k+1}}$.

  ($\ce$) Now we show that \eqref{eqn:use_of_tileset} is a subset of
  \eqref{eqn:use_pres_of_inter}. This direction is easier.

    Let $x \in \hull_k \less \hull_{k+1}$ so that $x=\simt_j(y)$
    for some $y \in \hull_{k-1}$. We know $y \notin \hull_k$, because
    otherwise $x = \simt_j(y) \in \simt_j(\hull_k) \ci \hull_{k+1}$.
    Thus $y \in \hull_{k-1} \less \hull_k$, which implies $x \in
    \simt(\hull_{k-1} \less \hull_k)$.
    This establishes $\hull_k \less \hull_{k+1} \ci \simt(\hull_{k-1}
    \less \hull_k)$; taking closures completes the proof of the equality
    \eqref{eqn:propagation_of_tilesets}.
  \end{proof}
\end{theorem}

Consider the set $\hull_k \less \hull_{k+1}$ and the set
$\simt^k(\bigcup \gen_q)$. The following result states that, modulo
some boundary points, these two sets are the same. In terms of
Figure~\ref{fig:koch-contractions}, for example, this means that the
difference between two successive stages in the left column is
roughly equal to the corresponding stage in the right column.

\begin{theorem}
  \label{thm:propagation_of_tiles}
  The tilesets can be recovered as the closure of the images of the
  generators under the action of \simt, that is,
  \begin{equation}
    \label{eqn:tileset1 as union}
    \tileset_{k} = \cj{\bigsqcup\nolimits_{q=1}^Q \simt^{k-1}(\gen_q)}.
  \end{equation}
  \begin{proof}
    First, observe that Corollary~\ref{cor:hull-less-simt(hull)-is-nice}
    gives
    \begin{align}
      \label{eqn:tileset1_is_union_of_closed_gens}
      \bigcup_{q=1}^Q \cj{\gen_q}
      = \cj{\bigsqcup\nolimits_{q=1}^Q \gen_q}
      = \cj{\inter(\hull \less \simt(\hull))}
      = \cj{\hull \less \simt(\hull)}
      = \tileset_1.
    \end{align}
    Now take $\simt^{k-1}$ of both sides, using
    Lemma~\ref{lemma:simt_preserves_closures} on the left and
    Theorem~\ref{thm:propagation_of_tilesets} on the right, to
    obtain the conclusion:
    \begin{align}
      \label{eqn:tilesetk+1_is_simt_of_union_of_closed_gens}
      \cj{\bigsqcup\nolimits_{q=1}^Q \simt^{k-1}\left(\gen_q\right)}
      = \simt^{k-1}\left(\cj{\bigsqcup\nolimits_{q=1}^Q \gen_q}\right)
      = \simt^{k-1}(\tileset_1)
      = \tileset_{k}.
    \end{align}
    The union at left is disjoint \bc each $\simt_j$ is \inj, $\gen_q
    \ci \inter \hull$, and the tileset condition
    \eqref{eqn:tileset_condition} prohibits overlaps of interiors.
  \end{proof}
\end{theorem}

\begin{theorem}
  \label{thm:the tiling}
  The collection $\tiling = \{\simt_w(\gen_q)\}$ is a tiling of
  $\hull \less \attr$ by open sets:
  \begin{equation}
    \label{eqn:hull_is_a_tiling}
    \hull = \cj{\bigcup \tile_n} = \cj{\bigcup \simt_w(\gen_q)}.
  \end{equation}
  \begin{proof}
    (i) To see the forward inclusion of \eqref{eqn:hull_is_a_tiling},
      let $x \in \hull$. If $x \notin \attr$, we can find $k$
      \st $x \in \hull_{k-1} \less \hull_k \ci  \tileset_k$, by
      Cor.~\ref{cor:hulls_converge_to_F}.
      By Thm.~\ref{thm:propagation_of_tiles}, it follows that $x \in
      \cj{\bigsqcup\nolimits_{q=1}^Q \simt^{k-1}(\gen_q)}$ and we are
      done. Now suppose that $x \in \attr$, and let $B_i$ be the open
      ball around $x$ of radius $1/i$.
      By Cor.~\ref{cor:hulls_converge_to_F}
      again, we can find $x_i \in B_i \cap (\hull \less \attr)$.
      The previous argument shows $x_i \in \cj{\bigcup
      \simt_w(\gen_q)}$, and hence the same holds for $x = \lim
      x_i$.
      The reverse inclusion is obvious from
      Theorem~\ref{thm:successor_is_always_contained} and the
      definition of the tiles as subsets of the $\hull_k$, in
      \eqref{def:generators}.

    (ii) To see that the tiles are disjoint, note first that the
      generators are disjoint by definition.
      Suppose $\tile_n$ and $\tile_m$ are both in the same
      tileset $\tileset_k$. Then
      \eqref{eqn:tilesetk+1_is_simt_of_union_of_closed_gens}
      shows that they are disjoint.
      Now suppose $\tile_n \ci \tileset_k$ and $\tile_m \ci
      \tileset_\ell$, where $k < \ell$. Then $\tile_n$ is disjoint from
      $\hull_k$ by definition of $\tileset_k$, and it follows from
      Theorem~\ref{thm:successor_is_always_contained} that
      $\tile_n$ is disjoint from $\hull_\ell$ for all $\ell \geq k$.
      (See, e.g., Figure \ref{fig:koch-contractions}.)

      It is also clear that $\tile_n \cap \hull_k = \es$ implies
      that $\tile_n \cap \attr = \es$, so no tiles intersect the
      attractor \attr. Thus, \tiling is an open tiling of $\hull
      \less \attr$.
  \end{proof}
\end{theorem}

\begin{cor}
  \label{cor:tiling_is_subselfaffine}
  The tiling \tiling is subselfaffine in that $\simt(\tiling) = \tiling \less \bigsqcup_q \gen_q$.
\end{cor}

\begin{remark}
  \label{ssec:Properties of the generators}
  What kinds of generators are possible? In general, this is a
  difficult question to answer; it is explored in detail in
  \cite{TFG}. The generators inherit many geometric properties from
  the convex hull $\hull = [\attr]$ and may therefore have a finite or
  infinite number of nonregular boundary points. In fact, by an
  observation of \cite{GAST1}, it is possible (even generic) for the
  boundary of a 2-\dimnl generator to be a piecewise $C^1$ curve.
  However, it is impossible for it to be a piecewise $C^2$ curve,
  unless it is polyhedral.
\end{remark}

\begin{remark}
  \label{rem:why_the_convex_hull}
  One might ask why the convex hull plays such a unique role in the
  construction of the tiling. There may exist other sets which are
  suitable for initiating the construction; however, some properties
  seem to make the convex hull the natural choice:
  \begin{enumerate}
    \item Any convex set (which is not a singleton set) is the closure
      of its relative interior 
      (as shown in
      the proof of Cor.~\ref{cor:Cks_are_nondegenerate}).
    \item The affine image of a convex set is convex. Consequently,
      $\simt(C) \ci C$ as in Theorem~\ref{thm:successor_is_always_contained}.
    \item The convex hull of \attr obviously contains \attr.
  \end{enumerate}
  Note that nonlinear maps do not preserve convexity, and so the
  convex hull would likely not be appropriate for constructing a
  tiling in the case that the mappings are not affine.
  
  Bandt et al. introduce the notion of \emph{central open set} in \cite{Bandt:OnOSC}. This provides a more natural (but less intuitive) feasible open set which is nonempty precisely when the open set condition is satisfied, and its closure may provide a substitute for the convex hull of the attractor. This possibility is considered in \cite{GeometryOfSST}.
\end{remark}

\pgap

\emph{Acknowledgements.} I am grateful to my advisor, Michel
Lapidus, for his advice, encouragement, and many fruitful
discussions. I would also like to thank Bob Strichartz and Steffen
Winter for comments on a preliminary version of the paper, and
several excellent references. Through many detailed and very helpful
remarks, the referee contributed considerably to the improvement of
this paper.

The idea for the tiling was inspired by the approach of Lapidus and
van Frankenhuijsen \cite[Ch.~2]{FGCD} in the 1-\dimnl case, and also
partially by trying to find a covering reminiscent to that of
Whitney, as in \cite{Stein:SingInt}, but naturally suited to \simt.
Michel Lapidus suggested that I investigate Whitney coverings.

\pgap

\bibliographystyle{harvard}

\pgap

\par

\end{document}